\documentclass{amsart}
\usepackage{graphicx}
\theoremstyle{definition}

\theoremstyle{remark}

\numberwithin{equation}{section}

\begin{document}

\title[]{Differential Algebras with Dense Singularities on Manifolds}%

\author{Elem\'{e}r E Rosinger}%

\address{Department of Mathematics and Applied Mathematics, University of Pretoria, 0002,
Pretoria, South Africa}%

\email{eerosinger@hotmail.com}

% ----------------------------------------------------------------
\begin{abstract}

Recently the space-time foam differential algebras of generalized functions with \textit{dense}
singularities were introduced, motivated by the so called space-time foam structures in
General Relativity with dense singularities, and by Quantum Gravity. A variety of applications
of these algebras has been presented, among them, a global Cauchy-Kovalevskaia theorem, de
Rham cohomology in abstract differential geometry, and so on. So far the space-time foam
algebras have only been constructed on Euclidean spaces. In this paper, owing to their
relevance in General Relativity among others, the construction of these algebras is extended
to arbitrary finite dimensional smooth manifolds. Since these algebras contain the Schwartz
distributions, the extension of their construction to manifolds also solves the long
outstanding problem of defining distributions on manifolds, and doing so in ways compatible
with nonlinear operations. Earlier, similar attempts were made in the literature with respect
to the extension of the Colombeau algebras to manifolds, algebras which also contain the
distributions. These attempts have encountered significant technical difficulties, owing to
the growth condition type limitations the elements of Colombeau algebras have to satisfy near
singularities. Since in this paper no any type of such or other growth conditions are required
in the construction of space-time foam algebras, their extension to manifolds proceeds in a
surprisingly easy and natural way. It is also shown that the space-time foam algebras form a
fine and flabby sheaf, properties which are important in securing a considerably large class
of singularities which generalized functions can handle.

\vspace{1cm}

\begin{quote}

"We do not possess any method at all to derive systematically solutions that are free of
singularities..."

\medskip

\hspace{4cm} Albert Einstein

\hspace{4.05cm}\textit{The Meaning of Relativity}

\hspace{4.1cm}Princeton Univ. Press, 1956, p. 165

\end{quote}

\end{abstract}

\maketitle

\newpage

% ----------------------------------------------------------------
\section{\bf Review of Space-Time Foam Algebras}

\subsection{Introductory Remarks}

There are \textit{five} features of interest in the construction in this paper of the
space-time foam differential algebras of generalized functions on arbitrary finite dimensional
smooth manifolds, namely

\begin{itemize}

\item the motivation for, and application of such space-time foam algebras is in General
Relativity and Quantum Gravity, see Mallios \& Rosinger [1-3], Mallios [1,2],

\item these algebras can handle dense singularities, and in fact, the cardinal of the
singularities can be larger than than of the regular or nonsingular points,

\item in the neighbourhood of singularities, there are no growth or other type of restrictions
on the generalized functions in space-time foam algebras,

\item unlike previous recent constructions of generalized functions or distributions on
manifolds, Grosser et.al., the construction in this paper of the space-time foam algebras on
such manifolds turns out to be surprisingly easy and natural,

\item unlike various spaces of usual or generalized functions, among the latter, the Schwartz
distributions, or for that matter, the Colombeau algebras, the space-time foam algebras in
this paper can form a flabby sheaf, a property which is closely related to their ability to
handle very large classes of singularities, see subsection \ref{subsec2.3}., or for more
details, Rosinger [7, pp. 173-187].

\end{itemize}

In Rosinger [8-10], motivated among others by the so called space-time foam structures with
dense singularities in General Relativity the class of space-time foam differential algebras
of generalized functions was constructed. A major interest in large, possibly dense sets of
singularities has a significant literature, as for instance in Finkelstein, Geroch [1,2],
Heller [1-3], Heller \& Sasin [1-3], Heller \& Multarzynski \& Sasin, Gurszczak \& Heller, or
Mallios [1-6], Mallios \& Rosinger [1-3], as well as the earlier literature cited there. \\
We note in this regard that according to the strongest earlier results, see Heller [2], Heller
\& Sasin [2], one could only deal with a singularity given by one single closed nowhere dense
set, which in addition had to be on the boundary of the domain. \\
In Mallios \& Rosinger [1] which, except for Rosinger [10] and Mallios \& Rosinger [2,3], does
so far treat the most general type of singularities, the family of singularities could already
contain all closed nowhere dense subsets, and they could be situated anywhere, and not only on
the boundary. Finally, in Rosinger [10] and Mallios \& Rosinger [2,3], the largest class of
singularities so far, namely such as in this paper, thus in particular, dense singularities as
well, were treated.

\bigskip

For earlier developments regarding the possible treatment of singularities in a differential
geometric context one can consult for instance, Sikorsky, Kirillov [1,2], Mostow or Souriau
[1,2]. And it should be mentioned that, as seen in Finkelstein, and especially in Geroch [1,2],
the issue of singularities has for a longer time been of fundamental importance in General
Relativity.

\bigskip

A recent significant effort in this regard can be found in Grosser et.al., where the Colombeau
differential algebras of generalized functions are alone employed. As it happens, however, the
considerable extent of the technical complications involved in the approach in Grosser et.al.
when defining on manifolds the Colombeau generalized functions is in clear contradistinction
with the natural ease in this paper when the space-time foam algebras with dense - therefor,
far worse singularities - are defined on manifolds.

\medskip

Indeed, one of the main results in this paper is what in fact amounts to a surprising and so
far unprecedented natural ease in defining on arbitrary finite dimensional smooth manifolds
differential algebras of generalized functions with such large classes of singularities as
those of the space-time foam algebras, singularities which can not only be dense in the
respective manifolds, but can also have a cardinality larger than that of the regular or
nonsingular points.

\medskip

A further interest in the space-time foam algebras comes from recent research in Quantum
Gravity, see Mallios [2], where such algebras prove to play a fundamental role.

\bigskip

The class of space-time foam algebras proves to be another specific instance of the
differential algebras of generalized functions defined and characterized earlier in Rosinger
[1-7], see also subsection \ref{subsec1.7}. for certain details on that characterization. And
the main interest of this class is in the fact that, this time, the \textit{singularity} sets
of the generalized functions which are elements of these algebras can be \textit{arbitrary},
provided that their complementaries, that is, the set of nonsingular or regular points, is
still dense in the domain where the respective generalized functions are defined.
\\
For instance, if the domain is the real line $\mathbb{R}$, the singularity set can be that of
all irrational numbers, since its complementary, the set of rational numbers, happens still to
be dense in $\mathbb{R}$. In this way, the set of singularities of a generalized function can
itself be \textit{dense} and \textit{uncountable}, while the set of its regular points only
must be dense, even if only countable. It follows, among others, that the cardinal of the set
of singularities can be \textit{larger} than that of the set of regular points.

\bigskip

These space-time foam algebras are thus far more powerful in handling singularities than the
most powerful earlier such algebras $A_{nd}(X)$, called the nowhere dense differential
algebras of generalized functions, which were introduced and used in Rosinger [3-7,10], see
also Mallios \& Rosinger [1].  Indeed, these latter algebras could only handle singularities
concentrated on closed and nowhere dense sets. And such sets - although may have arbitrary
large positive Lebesgue measures, see Oxtoby - are nevertheless far from being dense in the
domains of the respective generalized functions.
\\
However, in spite of such limitations on the singularities handled by the nowhere dense
algebras, one could already obtain in them \textit{global} versions of the Cauchy-Kovalevskaia
theorem, see Rosinger [4-7,10]. So far, this result could not be replicated in any other of
the earlier algebras of generalized functions, including those of Colombeau, which themselves
enter as a particular case of the mentioned characterization in Rosinger [1-7], and as a
consequence, are a particular case of the algebras constructed in Rosinger [1-10], see also in
this regard Grosser et.al. [p. 7].

\bigskip

The above may suggest the extent of the increased power of the space-time foam algebras in
handling singularities. And as was expected, in Rosinger [10] it was shown that these algebras
can also deliver global existence of solution results in the classical Cauchy-Kovalevskaia
setup. Moreover, the global generalized solutions thus obtained will have \textit{better}
regularity properties than those obtained earlier in the nowhere dense algebras.

\subsection{Construction of Space-Time Foam Algebras}\label{subsec1.2}

Before going to the construction on arbitrary finite dimensional smooth manifolds, we consider,
as preparation, the case of nonvoid open subsets in usual Euclidean spaces. Let us therefore
have our domain for generalized functions be given by any non-void open subset $X$ of
$\mathbb{R}^n$. We shall consider various families of singularities in $X$, each such family
being given by a corresponding set $\mathcal{S}$ of subsets $\Sigma \subset X$, with each such
subset $\Sigma$ describing a possible set of singularities of a certain given generalized
function.

\medskip

The \textit{largest} family of singularities $\Sigma\subset X$ which we can consider is given
by

\begin{equation}\label{label1}
\mathcal{S}_{\mathcal{D}}(X)=\{~ \Sigma\subset X ~|~ X\setminus\Sigma
                          \textnormal{ is dense in } X ~\}
\end{equation}

The various families $\mathcal{S}$ of singularities $\Sigma\subset X$ which we shall deal with,
will therefore each satisfy the condition $\mathcal{S}\subseteq \mathcal{S}_{\mathcal{D}}(X)$.

\medskip

In this paper, as in Rosinger [8-10] and Mallios \& Rosinger [2,3], the family $\mathcal{S}$
of singularities can be any subset of $\mathcal{S}_{\mathcal{D}}(X)$ in (\ref{label1}). Among
other ones, two such families which will be of interest are the following

\begin{equation}\label{label2}
\mathcal{S}_{nd}(X)=\{~\Sigma \subset X ~|~ \Sigma
\textnormal{ is closed and nowhere dense in } X ~\}
\end{equation}

and

\begin{equation}\label{label3}
\mathcal{S}_{\textit{Baire I}}(X)=\{~\Sigma \subset X ~|~
\Sigma \textnormal{ is of first Baire category in } X ~\}
\end{equation}

Obviously

\begin{equation}\label{label4}
\mathcal{S}_{nd}(X) ~\subset~ \mathcal{S}_{\textit{Baire I}}(X) ~\subset~
\mathcal{S}_{\mathcal{D}}(X)
\end{equation}

\subsection{Asymptotically Vanishing Ideals}

The construction of the space-time foam algebras, first introduced in Rosinger [8], has
\textit{two} basic ingredients involved.

\medskip

First, we take any family $\mathcal{S}$ of singularity sets $\Sigma\subset X$, family which
satisfies the following two conditions

\begin{equation}\label{label5}
  \begin{array}{l}
        \forall \quad \Sigma \in \mathcal{S}~: \\ \\
        \quad \quad X \setminus \Sigma \textnormal{ is dense in } X
  \end{array}
\end{equation}

and

\begin{equation}\label{label6}
  \begin{array}{l}
        \forall \quad \Sigma,~ \Sigma ' \in \mathcal{S}~: \\ \\
        \exists \quad \Sigma '' \in \mathcal{S}~: \\ \\
        \quad \quad \Sigma \cup \Sigma ' \subseteq \Sigma ''
  \end{array}
\end{equation}

\bigskip

Clearly, we shall have the inclusion $\mathcal{S}\subseteq \mathcal{S}_{\mathcal{D}}(X)$ for
any such family $\mathcal{S}$. Also, it is easy to see that both families $\mathcal{S}_{nd}
(X)$ and $\mathcal{S}_{\textit{Baire I}}(X)$ satisfy conditions (\ref{label5}) and
(\ref{label6}).

\medskip

Now, as the second ingredient, and so far independently of any $\mathcal{S}$ above, we take
any right directed partial order $L=(\Lambda,\leq)$. In other words, $L$ is such that for each
$\lambda,\, \lambda ' \in \Lambda$ there exists $\lambda ''\in \Lambda$, for which $\lambda,
\lambda ' \leq \lambda ''$. Here we note that the choice of $L$ may at first appear to be
completely independent of $\mathcal{S}$, yet in certain specific instances the two may be
somewhat related, with the effect that $\Lambda$ may have to be large, see Rosinger [10].

\bigskip

Although we shall only be interested in singularity sets $\Sigma \in \mathcal{S}_{\mathcal{D}}
(X)$, the following \textit{ideal} can in fact be defined for any $\Sigma \subseteq X$. Indeed,
let us denote by

\begin{equation}\label{label7}
\mathcal{J}_{L,~\Sigma}(X)
\end{equation}

the \textit{ideal} in $(\mathcal{C}^{\infty}(X))^{\Lambda}$ of all the sequences of smooth
functions indexed by $\lambda \in \Lambda$, namely, $w = (~ w_{\lambda} ~|~ \lambda\in
\Lambda ~)\in (\mathcal{C}^{\infty}(X))^{\Lambda}$, sequences which \textit{outside} of the
singularity set $\Sigma$ will satisfy the \textit{asymptotic vanishing} condition

\begin{equation}\label{label8}
      \begin{array}{l}
            \forall \quad x \in X \setminus \Sigma~: \\ \\
            \exists \quad \lambda \in \Lambda~: \\ \\
            \forall \quad \mu \in \Lambda,~ \mu \geq \lambda~: \\ \\
            \forall \quad p \in \mathbb{N}^n~: \\ \\
            \quad \quad D^p w_{\mu}(x)=0
       \end{array}
\end{equation}

This means that the sequences of smooth functions $w = (~ w_\lambda ~|~ \lambda \in \Lambda
~)$ in the ideal $\mathcal{J}_{L,~\Sigma}(X)$ will in a way \textit{cover} with their support
the singularity set $\Sigma$, and at the same time, they vanish outside of it asymptotically,
together with all their partial derivatives.
\\
In this way, the ideal $\mathcal{J}_{L,~\Sigma}(X)$ carries in an \textit{algebraic} manner
the information on the singularity set $\Sigma$.  Therefore, a \textit{quotient} space in
which the factorization is made with such an ideal may in certain ways \textit{do away with}
singularities, and do so through purely algebraic means, see (\ref{label11}), (\ref{label12})
below.

\bigskip

We note that the assumption about $L=(\Lambda,\leq)$ being right directed is used in proving
that $\mathcal{J}_{L,~\Sigma}(X)$ is indeed an ideal, more precisely that, for $w, w' \in
\mathcal{J}_{L,~\Sigma}(X)$, we have $w+w' \in \mathcal{J}_{L,~\Sigma}(X)$.

\medskip

Now, it is easy to see that for $\Sigma,~ \Sigma\, ' \subseteq X$ we have

\begin{equation}\label{label9}
\Sigma \subseteq \Sigma\, ' ~~\Longrightarrow~~ \mathcal{J}_{L,~\Sigma}(X)
\subseteq \mathcal{J}_{L,~\Sigma\, '}(X)
\end{equation}

In this way, for any family $\mathcal{S}$ of singularity sets $\Sigma\subset X$ satisfying
(\ref{label6}), it follows that

\begin{equation}\label{label10}
\mathcal{J}_{L,~\mathcal{S}}(X) ~=~ \bigcup_{\Sigma \in ~\mathcal{S}}~
                        \mathcal{J}_{L,~\Sigma}(X)
\end{equation}

is also an \textit{ideal} in $(\mathcal{C}^{\infty}(X))^{\Lambda}$

\bigskip

It is important to note that for suitable choices of the right directed partial orders $L$,
the ideals $\mathcal{J}_{L,\Sigma}(X)$, with $\Sigma \in \mathcal{S}_{\mathcal{D}}(X)$, are
\textit{nontrivial}, that is, they do not reduce to the zero ideal $\{~0~\}$. Thus in view of
(\ref{label10}), the same will hold for the ideals $\mathcal{J}_{L,~\mathcal{S}}(X)$. In
Rosinger [10, section 2] further details are presented in the case of general singularity sets
$\Sigma \in \mathcal{S}_{\mathcal{D}}(X)$ , when the respective right directed partial orders
$L$, which give the nontriviality of the ideals $\mathcal{J}_{L,~\Sigma}(X)$, prove to be
rather large, and in particular, uncountable. In subsection 1.9. we shall show that in the
case of the singularity sets $\Sigma$ in
$\mathcal{S}_{nd}(X)$ and $\mathcal{S}_{\delta \textit{ Baire I}}(X)$ - with the latter, see
(\ref{label28}) below, a suitable subset of $\mathcal{S}_{\textit{Baire I}}(X)$ - one can have
the ideals $\mathcal{J}_{L,~\Sigma}(X)$ nontrivial even for $L = \mathbb{N}$, respectively, $L
= \mathbb{N} \times \mathbb{N}$, that is, with $L$ still countable. On the other hand, in view
of (\ref{label13}), (\ref{label15}) below, the mentioned ideals cannot become too large either.

\subsection{Foam Algebras}

In view of the above, for $\Sigma \subseteq X$ we can define the \textit{quotient} algebra

\begin{equation}\label{label11}
B_{L,~\Sigma}(X) = (\mathcal{C}^{\infty}(X))^{\Lambda} /
\mathcal{J}_{L,~\Sigma}(X)
\end{equation}

However, we shall only be interested in singularity sets $\Sigma \in \mathcal{S}_{\mathcal{D}}
(X)$, that is, for which $X \setminus \Sigma$ is dense in $X$, and in such a case the
corresponding algebra $B_{L,~\Sigma}(X)$ will be called a \textit{foam algebra}.

\subsection{Multi-foam Algebras}

With any family $\mathcal{S}$ of singularities satisfying (\ref{label5}) and (\ref{label6}),
we can, based on (\ref{label10}), associate the \textit{multi-foam algebra}

\begin{equation}\label{label12}
B_{L,~\mathcal{S}}(X)=(\mathcal{C}^{\infty}(X))^{\Lambda}/
\mathcal{J}_{L,~\mathcal{S}}(X)
\end{equation}

\subsection{Space-Time Foam Algebras}

The foam algebras and the multi-foam algebras introduced above will for the sake of simplicity
be called together \textit{space-time foam algebras.}

\medskip

Clearly, if the family $\mathcal{S}$ of singularities consists of one single singularity set
$\Sigma \in \mathcal{S}_{\mathcal{D}}(X)$, that is, $\mathcal{S} = \{~ \Sigma ~\}$, then
conditions (\ref{label5}), (\ref{label6}) are trivially satisfied, and in this particular case
the concepts of foam and multi-foam algebras are identical, in other words, $B_{L, \{\Sigma\}}
(X)=B_{L, \Sigma}(X)$. This means that the concept of multi-foam algebra is more general than
that of foam algebra.

\bigskip

It is clear from their quotient construction that the space-time foam a1gebras are associative
and commutative. Later, we shall see that they are as well unital.

\medskip

However, the above constructions can easily be extended to the case when, instead of real
valued smooth functions, we use smooth functions with values in an arbitrary \textit{normed
algebra.} In such a case the resulting space-time foam algebras will still be associative, but
in general they may be \textit{non-commutative}.

\subsection{Space-Time Foam Algebras as Differential Algebras of Generalized Functions}\label{subsec1.7}

The reason why we restrict ourself to singularity sets $\Sigma \in \mathcal{S}_{\mathcal{D}}
(X)$, that is, to subsets $\Sigma \subset X$ for which $X \setminus \Sigma$ is dense in $X$,
is due to the implication

\begin{equation}\label{label13}
X \setminus \Sigma \textnormal{ is dense in } X ~~\Longrightarrow~~ \mathcal{J}_{L,~\Sigma}(X)
\cap \mathcal{U}_{\Lambda}^{\,\infty}(X)~=~ \{~ 0 ~\}
\end{equation}

where $\mathcal{U}_{\Lambda}^{\,\infty}(X)$ denotes the \textit{diagonal} of the power
$(\mathcal{C}^{\infty}(X))^{\Lambda}$, namely, it is the set of all $u(\psi)=(~ \psi_{\lambda}
~|~ \lambda \in \Lambda ~)$, where $\psi_{\lambda}=\psi$, for $\lambda \in \Lambda$, while
$\psi$ ranges over $\mathcal{C}^{\infty}(X)$. In this way, we have the algebra isomorphism

$$\mathcal{C}^{\infty}(X) \ni \psi \longmapsto u(\psi) \in
                              \mathcal{U}^{\,\infty}_{\Lambda}(X)$$

\bigskip

This implication (\ref{label13}) follows immediately from the asymptotic vanishing condition
(\ref{label8}), and it means that the ideal $\mathcal{J}_{L,~\Sigma}(X)$ is \textit{off
diagonal.} And in this implication the fact that $X \setminus \Sigma \textnormal{ is dense
in } X$ is essential. This is how we arrive at the restriction on singularity sets $\Sigma$
that their complementary $X \setminus \Sigma$ must be dense in $X$. However, as seen, this is
such a mild restriction that it allows the singularity sets $\Sigma$ to be themselves \textit
{dense} in $X$, and in fact, to have a cardinal \textit{larger} than that of the regular or
nonsingular points.

\medskip

The importance of (\ref{label13}) is that, for $\Sigma \in \mathcal{S}_{\mathcal{D}}(X)$, it
gives the following \textit{algebra embedding} of the smooth functions into foam algebras

\begin{equation}\label{label14}
\mathcal{C}^{\infty}(X) \ni \psi \longmapsto u(\psi) + \mathcal{J}_{L,~\Sigma}(X) \in
B_{L,~\Sigma}(X)
\end{equation}

Now in view of (\ref{label10}), it is easy to see that (\ref{label13}) will yield the \textit
{off diagonality} property as well

\begin{equation}\label{label15}
\mathcal{J}_{L,~\mathcal{S}}(X) \cap
\mathcal{U}^{\,\infty}_{\Lambda}(X)=\{~ 0 ~\}
\end{equation}

and thus similar with (\ref{label14}), we obtain the \textit{algebra embedding} of smooth
functions into multi-foam algebras

\begin{equation}\label{label16}
\mathcal{C}^{\infty}(X) \ni \psi \longmapsto u(\psi) + \mathcal{J}_{L,~\mathcal{S}}(X) \in
B_{L,~\mathcal{S}}(X)
\end{equation}

The algebra embeddings (\ref{label14}), (\ref{label16}) mean that the foam and multi-foam
algebras are in fact \textit{algebras of generalized functions.} Also they mean that the foam
and multi-foam algebras are unital, with the respective unit elements $u(1) + \mathcal{J}_{L,
~\Sigma}(X$) and $u(1) + \mathcal{J}_{L,~\mathcal{S}}(X)$.

\bigskip

Further, the asymptotic vanishing condition (\ref{label8}) also implies quite obviously that,
for $\Sigma \subseteq X$  we have

\begin{equation}\label{label17}
D^p \,~\mathcal{J}_{L, \Sigma}(X) \subseteq \mathcal{J}_{L,~\Sigma}(X), \quad
\textnormal{ for } p\in \mathbb{N}^n
\end{equation}

where $D^p$ denotes the termwise \textit{p}-th order partial derivation of sequences of smooth
functions, applied to each such sequence in the ideal $\mathcal{J}_{L,~\Sigma}(X)$. Thus in
view of (\ref{label10}), we obtain

\begin{equation}\label{label18}
D^p \, \mathcal{J}_{L,~\mathcal{S}}(X)\subseteq \mathcal{J}_{L,~\mathcal{S}}(X), \quad
\textnormal{ for } p\in \mathbb{N}^n
\end{equation}

Now (\ref{label17}), (\ref{label18}) mean that the the foam and multi-foam algebras are in
fact \textit{differential algebras}, namely

\begin{equation}\label{label19}
D^p \, B_{L,~\Sigma}(X) \subseteq B_{L,~\Sigma}(X), \quad
\textnormal{ for } p\in \mathbb{N}^n
\end{equation}

where $\Sigma \in \mathcal{S}_{\mathcal{D}}(X)$, while we also have

\begin{equation}\label{label20}
D^p \, B_{L,~\mathcal{S}}(X)\subseteq B_{L,~\mathcal{S}}(X), \quad
\textnormal{ for } p\in \mathbb{N}^n
\end{equation}

for $\mathcal{S} \subseteq \mathcal{S}_{\mathcal{D}}(X)$ which satisfy (\ref{label5}) and
(\ref{label6}).

\medskip

In this way we obtain that the foam and multi-foam algebras are \textit{differential algebras
of generalized functions.}

\medskip

Also, the multi-foam algebras contain the Schwartz distributions, that is, we have the
\textit{linear embeddings} which respect the arbitrary partial derivation of smooth functions

\begin{equation}\label{label21}
\mathcal{D}\, ' (X) ~\subset~ B_{L,~\Sigma}(X), \quad \textnormal{ for }
\Sigma \in \mathcal{S}_{\mathcal{D}}(X)
\end{equation}

\begin{equation}\label{label22}
\mathcal{D}\, ' (X) ~\subset~ B_{L,~\mathcal{S}}(X)
\end{equation}

Indeed, we mentioned the wide ranging purely algebraic characterization of all those quotient
type algebras of generalized functions in which one can embed linearly the Schwartz
distributions, a characterization first given in 1980, see Rosinger [4, pp. 75-88], Rosinger
[5, pp. 306-315], Rosinger [6, pp. 234-244].  According to that characterization - which also
contains the Colombeau algebras as a particular case - the necessary and sufficient condition
for the existence of the linear embedding (\ref{label21}) is precisely the off diagonality
condition in (\ref{label13}). Similarly, the necessary and sufficient condition for the
existence of the linear embedding (\ref{label22}) is precisely the off diagonality condition
(\ref{label15}).

\bigskip

One more property of the foam and multi-foam algebras will prove to be useful.  Namely, in
view of (\ref{label10}), it is clear that, for every $\Sigma \in \mathcal{S}$, we have the
inclusion $\mathcal{J}_{L,~\Sigma}(X) \subseteq \mathcal{J}_{L,~\mathcal{S}} (X)$, and thus we
obtain the \textit{surjective algebra homomorphism}

\begin{equation}\label{label23}
B_{L,~\Sigma}(X) \ni w + \mathcal{J}_{L,~\Sigma}(X) \longmapsto w +
\mathcal{J}_{L,~\mathcal{S}}(X) \in B_{L,~\mathcal{S}}(X)
\end{equation}

As we shall see in the next subsection, (\ref{label23}) can naturally be interpreted as
meaning that the typical generalized functions in $B_{L,~\mathcal{S}}(X)$ are \textit{more
regular} than those in $B_{L,~\Sigma}(X)$.

\subsection{Regularity of generalized Functions}\label{subsec1.8}

One natural way to interpret (\ref{label23}) in the context of generalized functions is the
following. Given two spaces of generalized functions $E$ and $F$, such as for instance

\begin{equation}\label{label24}
\mathcal{C}^{\infty}(X) ~\subset~ E ~\subset~ F
\end{equation}

then the \textit{larger} the space $F$, the \textit{less regular} its typical element can
appear to be, when compared with those of $E$. By the same token, the \textit{smaller} the
space $E$, the \textit{more regular} one can consider its typical elements, when compared to
those of $F$. Similarly, given a \textit{surjective} mapping

\begin{equation}\label{label25}
E ~~\longrightarrow~~ F
\end{equation}

one can again consider that the typical elements of $F$ are at least as \textit{regular} as
those of $E$ .

\medskip

In this way, in view of (\ref{label23}), we can consider that, owing to the given \textit
{surjective} algebra homomorphism, the typical elements of the multi-foam algebra $B_{L,~
\mathcal{S}}(X)$ can be seen as being \textit{more regular} than the typical elements of the
foam algebra $B_{L,~\Sigma}(X)$.

\bigskip

Furthermore, the algebra $B_{L,~\mathcal{S}}(X)$ is obtained by factoring the same
$(\mathcal{C}^{\infty}(X))^{\Lambda}$ as in the case of the algebra $B_{L,~\Sigma}(X)$, this
time however by the significantly \textit{larger} ideal $\mathcal{J}_{L,~\mathcal{S}}(X)$,
an ideal which, unlike any of the individual ideals $\mathcal{J}_{L,~\Sigma}(X)$, can
simultaneously deal with \textit{all} the singularity sets $\Sigma \in \mathcal{S}$, some, or
in fact, many of which singularity sets can be \textit{dense} in $X$. Further details related
to the connection between \textit{regularization} in the above sense, and on the other hand,
properties of \textit{stability, generality} and \textit{exactness} of generalized functions
and solutions can be found in Rosinger [4-6].

\subsection{ Nontriviality of Ideals}\label{subsec1.9}

Let us take any nonvoid singularity set $\Sigma \in \mathcal{S}_{nd}(X)$. Since $\Sigma$ is
closed, we can take a sequence of nonvoid open subsets $Y_l \subset X$, with $l \in
\mathbb{N}$, such that $\Sigma = \bigcap_{l\in \mathbb{N}} \, Y_l$. We can also assume that
the $Y_l$ are decreasing in $l$, since we can replace every $Y_l$ with the finite intersection
$\bigcap_{k\leq l}\,Y_k$. But for each $Y_l$, Kahn, there exists $\alpha_l \in \mathcal{C}^
{\infty}(X)$, such that $\alpha_l = 1$ on $\Sigma$, and $\alpha_l = 0$ on $X\setminus Y_l $.
Now in view of (\ref{label8}) it is easy to check that the resulting sequence of smooth
functions on $X$ satisfies

\begin{equation}\label{label26}
\alpha = (~\alpha_l ~|~ l\in \mathbb{N} ~)\in \mathcal{J}_{\mathbb{N},~\Sigma}(X)
\end{equation}

and clearly, $\alpha$ is not a trivial sequence, since $\phi \neq \Sigma \subseteq \textnormal
{supp } \alpha_l$, for $l \in \mathbb{N}$.

\bigskip

We can note, however, that the above argument leading to (\ref{label26}) need not necessarily
apply to subsets $\Sigma \subset X$ which are not closed, even if their closure is
nevertheless nowhere dense in $X$. In such a case one can use the more general method in
Rosinger [10, section 2], which will give nontrivial sequences similar to $\alpha$ above,
sequences whose index sets, however, will no longer be countable.

\bigskip

Let us take now any nonvoid singularity set $\Sigma \in \mathcal{S}_{\textit{Baire I}}(X)$.
Then there exists a sequence of closed and nowhere dense subsets $\Sigma_l \subset X$, with $l
\in \mathbb{N}$, such that

\begin{equation}\label{label27}
\Sigma \subseteq \bigcup_{l \in \mathbb{N}}\, \Sigma_l
\end{equation}

where the equality need not necessarily hold. Therefore, let us consider the subset of
$\mathcal{S}_{\textit{Baire I}}(X)$ denoted by

\begin{equation}\label{label28}
\mathcal{S}_{\delta\ \, \textit{Baire I}}(X)
\end{equation}

whose elements are all those singularity sets $\Sigma$ for which we have equality in
(\ref{label27}). Obviously, we can assume that the $\Sigma_l$ are increasing in $l$, since we
can replace each $\Sigma_l$ with the finite union $\bigcup_{k\leq l}\Sigma_k$.

\bigskip

Given now a nonvoid $\Sigma \in \mathcal{S}_{\delta \, \textit{Baire I}}(X)$ and a
corresponding representation $\Sigma = \bigcup_{l \in \mathbb{N}}\,\Sigma_l$, with suitable
closed and nowhere dense subsets $\Sigma_l \subseteq X$ which are increasing in $l$, we can,
as above, find for each $\Sigma_l$ a representation $\Sigma_l = \bigcap_{k\in \mathbb{N}} \,
Y_{l,\,k}$, with nonvoid open subsets $Y_{l,\,k} \subset X$.  Further, we can assume that for
$l, l', k, k' \in \mathbb{N}, \, l\leq l', \, k\leq k'$, we have $Y_{l,\,k} \supseteq Y_{l',\,
k'}$, since we can replace every $Y_{l,\,k}$ with the finite intersection $\bigcap_{l'\leq l,
\,k'\leq k}\, Y_{l',\,k'}$. And for each $Y_{l,\,k}$ we can find $\alpha_{l,\,k}\in \mathcal
{C}^{\infty}(X)$, such that $\alpha_{l,\,k}=1$ on $\Sigma_l$, while $\alpha_{l,\,k}=0$ on
$X\setminus Y_{lk}$.
\\
Let us now take $L=(\Lambda, \leq)$, where $\Lambda = \mathbb{N} \times \mathbb{N}$, while for
$(l,\,k),~ (l',\,k')\in \Lambda = \mathbb{N} \times \mathbb{N}$ we set $(l,\,k)\leq (l',\,k')$,
if and only if $l \leq l'$ and $k\leq k'$. Then (\ref{label8}) will easily give

\begin{equation}\label{label29}
\alpha= (~ \alpha_{l,\,k} ~|~ (l,k) \in \mathbb{N}\times \mathbb{N} ~) \in
\mathcal{J}_{\mathbb{N} \times \mathbb{N},~\Sigma}(X)
\end{equation}

And again, $\alpha$ is not a trivial sequence, since $\phi \neq \Sigma \subseteq \bigcup_{l
\in \mathbb{N}} \textnormal{ supp } \alpha_{l,\,k_l}$, for every given choice of $k_l \in
\mathbb{N}$, with $l \in \mathbb{N}$.

\bigskip

In case our singularity set $\Sigma$ belongs to $\mathcal{S}_{\textit{Baire I}}(X)$ but not to
$\mathcal{S}_{\delta \, \textit{Baire I}}(X)$, then the above approach need no longer work.
However, we can still apply the mentioned more general method in Rosinger [10, section 2] -
valid for every nonvoid singularity set $\Sigma\in \mathcal{S}_{\mathcal{D}}(X)$ - in order to
construct nontrivial sequences in $\mathcal{J}_{L,~\Sigma}(X)$, although this time the
corresponding index sets $\Lambda$ may be uncountable.

\subsection{Relations between Algebras with the Same Singularities}\label{subsec1.10}

The above, and especially subsection \ref{subsec1.8}., leads to the following question. Let us
assume given a certain nonvoid singularity set $\Sigma\in \mathcal{S}_{\mathcal{D}}(X)$. If we
now consider two right directed partial orders $L=(\Lambda, \leq)$ and $L' = (\Lambda ',
\leq)$, is there then any relevant relationship between the corresponding two foam algebras

\begin{equation}\label{label30}
B_{L,~\Sigma}(X) \textnormal{  and  } B_{L',~\Sigma}(X) \textnormal{~~?}
\end{equation}

A rather simple positive answer can be given in the following particular case. Let us assume
that $\Lambda$ is a \textit{cofinal} subset of $\Lambda '$, that is, the partial order on
$\Lambda$ is induced by that on $\Lambda '$, and in addition, we also have satisfied the
condition

\begin{equation}\label{label31}
     \begin{array}{l}
           \forall \quad \lambda ' \in \Lambda ~: \\ \\
           \exists \quad \lambda \in \Lambda ~: \\ \\
           \quad \quad \lambda ' \leq \lambda
     \end{array}
\end{equation}

Then considering the surjective algebra homomorphism

\begin{equation}\label{label32}
(\mathcal{C}^{\infty}(X))^{\Lambda '} \ni s' = (~ s'_{\lambda '} ~|~
\lambda ' \in \Lambda ' ~) ~~\stackrel{\rho}{\longmapsto}~~ s = (~ s'_{\lambda '} ~|~
\lambda ' \in \Lambda) \in (\mathcal{C}^{\infty}(X))^{\Lambda}
\end{equation}

and based on (\ref{label8}), one can easily note the property

\begin{equation}\label{label33}
\rho \, (~ \mathcal{J}_{\Lambda ' ,~\Sigma}(X) ~) \subseteq
\mathcal{J}_{\Lambda,~\Sigma}(X)
\end{equation}

In this way, one can obtain the \textit{surjective} algebra homomorphism of foam algebras,
given by

\begin{equation}\label{label34}
B_{\Lambda ',~\Sigma}(X) \ni s' + \mathcal{J}_{\Lambda ',~\Sigma}(X)
 ~~\stackrel{\rho}{\longmapsto}~~ \rho \, s' +\mathcal{J}_{\Lambda,~\Sigma}(X) \in
 B_{\Lambda,~\Sigma}(X)
\end{equation}

And in the terms of the interpretation in subsection \ref{subsec1.8}., the meaning of
(\ref{label34}) is that the foam algebra $B_{\Lambda,~\Sigma}(X)$ has its typical generalized
functions \textit{more} regular than those of $B_{\Lambda ',~\Sigma}(X)$. Thus in such terms,
foam algebras which correspond to a \textit{smaller cofinal} partial order can be seen as
\textit{less} regular.
\\
However, there may be many other kind of relationships between two partial orders $L$ and $L'$,
such as for instance in the case of $\mathbb{N}$ in (\ref{label26}) and $\mathbb{N} \times
\mathbb{N}$ in (\ref{label29}). Therefore the problem in (\ref{label30}) may in general
present certain difficulties.

\bigskip

Needless to say, similar results and comments hold in the case of the multi-foam, and thus
space-time foam algebras as well.

\section{ Algebras on Manifolds}

\subsection{ Extension to Manifolds}\label{subsec2.1}

We shall from now on assume that the domain $X$ of our generalized functions is no longer a
nonvoid open subset of $\mathbb{R}^n$, but an arbitrary n-dimensional smooth manifold, Kahn.
Also, we shall assume given an atlas $\mathcal{V}$ for the manifold $X$, atlas composed from
compatible charts $v: V \rightarrow \mathbb{R}^n$, with $V \in \mathcal{V}$ being open subsets
of $X$, the totality of which charts covers $X$, while $v$ are homeomorphisms.

\bigskip

What we can note immediately is that all of the definitions, constructions and properties in
the above subsections \ref{subsec1.2}. to \ref{subsec1.8}. - with the sole exception of
relations (\ref{label21}), (\ref{label22}) - will still remain valid, now that $X$ is no
longer a nonvoid open set of $\mathbb{R}^n$, but an arbitrary n-dimensional smooth manifold.
\\
Furthermore, the problem with relations (\ref{label21}), (\ref{label22}) is not with the
respective algebras, but with defining $\mathcal{D}\, ' (X)$, that is, with defining the
Schwartz distributions when $X$ is an arbitrary n-dimensional smooth manifold.

\medskip

Let us now return to the fact that all the mentioned constructions extend immediately to
arbitrary n-dimensional smooth manifolds $X$.
\\
In particular, this means that the nowhere dense differential algebras of generalized
functions $A_{nd}(X)$ can as well be extended in the above direct manner to the mentioned
manifolds $X$. Indeed, it is easy to see that if we take the right directed partial order
$L=(\Lambda, \leq)$ as given by the set of natural numbers $\mathbb{N}$, then in view of
(\ref{label4}) , (\ref{label12}), we obtain, see Rosinger [3-7,10]

\begin{equation}\label{label2.1}
A_{nd}(X) ~=~ B_{\,\mathbb{N},~\mathcal{S}_{nd}(X)}(X)
\end{equation}

In other words, the nowhere dense algebras can be seen as rather particular cases of the
space-time foam algebras.

\bigskip

Here it is important to point out the following.

\medskip

The above direct manner in which the space-time foam algebras can be defined on arbitrary
finite dimensional smooth manifolds has its reasons in the fact that the construction of these
algebras only involves two kinds of ingredients.  First, topology in Euclidean spaces, and in
ways which are actually of a local nature. Second, ring theoretic concepts related to
functions defined on manifolds.
\\
And as seen above, the properties involved in these kind of ingredients which are relevant in
the construction of the space-time foam algebras happen to be indifferent to the fact whether
they are considered on nonvoid open sets of Euclidean spaces, or on arbitrary finite
dimensional smooth manifolds.

\bigskip

On the other hand, the definition of the Schwartz distributions $\mathcal{D}'(X)$ on nonvoid
open subsets $X$ of Euclidean spaces involves essentially \textit{linear} concepts, be they
the duality of vector spaces, or weak convergence. And these concepts do not allow a direct
extension to arbitrary finite dimensional smooth manifolds.
\\
Also, in the case of the Colombeau algebras of generalized functions, their definition
essentially involves certain \textit{growth conditions.}  And as seen in a number of recent
attempts, among them Grosser et.al., these growth conditions create serious difficulties and
consequent technical complications when one tries to extend the Colombeau algebras to
arbitrary finite dimensional smooth manifolds.

\subsection{What About the Distributions ?}

As seen in (\ref{label22}), namely in the embedding

\begin{equation}\label{label2.2}
\mathcal{D}'(X) ~\subset~ B_{L,~\mathcal{S}}(X)
\end{equation}

which holds in the particular case of a  flat manifold $X$, that is, when $X$ is a nonvoid
open subset of $\mathbb{R}^n$, each of the space-time foam algebras $B_{L,\mathcal{S}}(X)$
contains the Schwartz distributions $\mathcal{D}'(X)$ as a linear vector subspace. Moreover,
the partial derivatives in $B_{L,\mathcal{S}}(X)$, see (\ref{label20}), coincide with the
usual ones, when restricted to smooth functions on $X$.

\bigskip

Now, as we have seen in subsection \ref{subsec2.1}., the extension to arbitrary finite
dimensional smooth manifolds $X$ of the right hand term in (\ref{label2.2}) is straightforward.
\\
On the other hand, as also mentioned in subsection \ref{subsec2.1}, the left hand term in
(\ref{label2.2}), due to its essentially linear nature, does not allow such a direct extension
to the same kind of manifolds.

\bigskip

In this way, in case $X$ is a finite dimensional smooth manifold, what one can save from
(\ref{label2.2}) is its \textit{local} version, that is

\begin{equation}\label{label2.3}
\mathcal{D}'(V) ~\subset~ B_{L,~\mathcal{S}}(X)|_{V} ~=~ B_{L,~\mathcal{S}|_{V}}(V)
\end{equation}

obtained by the \textit{restriction} of (\ref{label2.2}) to each chart $V\in \mathcal{V}$,
where we used the notation

\begin{equation}\label{label2.4}
\mathcal{S}|_{V} ~=~ \{~ \Sigma \cap V ~|~ \Sigma \in \mathcal{S} ~\}
\end{equation}

Here we should note that, as mentioned in Rosinger [7, pp. 5-8, 11-12, 173-187], one should
avoid rushing into a too early normative judgement about the way the long established
\textit{linear} theories of generalized functions - such as for instance that of the Schwartz
distributions - should relate to the still emerging, and far more complex and rich
\textit{nonlinear} theories.

\bigskip

In particular, \textit{two} aspects of such possible relationships still await a more
thoroughly motivated and clear settlement.
\\
First, the purely algebraic-differential type connections between Schwartz distributions and
the recent variety of differential algebras of generalized functions should be studied in more
detail. And since the main aim is to deal with generalized - hence, not smooth enough, but
rather singular - solutions, a main focus should be placed on the respective capabilities to
deal with singularities, see related comments in Rosinger [7, pp. 174-185]. In this regard let
us only mention the following. A fundamental property of various spaces of generalized
functions which is closely related to their capability to handle a large variety of
singularities is that such spaces should have a \textit{flabby} sheaf structure, see for
details Kaneko, or Rosingec [9,10]. However, the Schwartz distributions, the Colombeau algebra
of generalized functions, as well as scores of other frequently used spaces of classical or
generalized functions happen to fail being flabby sheaves. On the other hand, the nowhere
dense differential algebras of generalized functions, see (\ref{label2.1}) above, have a
structure of flabby sheaves, as shown for instance in Mallios \& Rosinger [1]. In subsection
\ref{subsec2.3}., we shall see that the same flabby sheaf structure also holds for the
space-time foam algebras.
\\
Needless to say, considerations which relate sheaf properties to generalized functions have so
far not been much prevalent in the main literature related to the latter, a somewhat earlier
remarkable exception in this regard being presented in Kaneko.
\\
However, with the recent emergence of Abstract Differential Geometry in Mallios [1,2], aimed
among others for Quantum Gravity, see also Mallios \& Rosinger [1-3], sheaf structures, more
precisely, sheaves of algebras, and specifically, sheaves of differential algebras of
generalized functions, among them the space-time foam algebras, turn out to play a fundamental
role.

\bigskip

The second aspect is possibly even more controversial. And it is so, especially because of the
historical phenomenon that the study of the \textit{linear} theories of generalized functions
has from its early modern stages in the 1930s been strongly connected with the then fast
emerging theories of linear topological structures.
\\
However, just as with the \textit{nonstandard} reals $^*\mathbb{R}$, so with the various
algebras of generalized functions, it appears that \textit{infinitesimal} type elements in
these algebras play an important role. And the effect is that if one introduces Hausdorff
topologies on these algebras, then, when these topologies are restricted to the more regular,
smooth, classical type functions, they inevitably lead to the trivial \textit{discrete}
topology on them, see related comments in the mentioned places in Rosinger [7].

\bigskip

In this way, due to the manifest lack of well enough researched and motivated arguments about
the variety of possible relationships between \textit{linear} and \textit{nonlinear} theories
of generalized functions, one should become aware of the question:

\bigskip

\begin{quote}

Are, indeed, relations like (\ref{label2.2}) - between the essentially linear distributions,
and on the other hand, various nonlinear generalized functions - the ones which one should try
to extend to manifolds ?

\end{quote}

\bigskip

Such a question has so far hardly been asked consciously, carefully and critically enough. Yet
it often tends to get a rather automatic affirmative answer.

\bigskip

Here we can also note that, compared with the space-time foam, or even the more particular
nowhere dense, differential algebras of generalized functions, the Colombeau algebras suffer
from several important limitations. Among them are the following two. First, the singularity
sets in the Colombeau algebras must be so \textit{small} as to be of zero Lebesgue measure,
and second, there are polynomial type \textit{growth conditions} which the \textit{generalized}
functions must satisfy in the neighbourhood of singularities.
\\
On the other hand, the earlier introduced nowhere dense algebras, Rosinger [3,4], as well as
the recent more general space time foam algebras, do not suffer from any of the above two
limitations. Indeed, the nowhere dense algebras allow singularities on arbitrary closed
nowhere dense sets, therefore, such singularity sets can have arbitrary large positive
Lebesgue measure, Oxtoby. Furthermore, the space-time foam algebras allow arbitrary
singularity sets, as long as the complementaries of such sets are dense. And in these algebras
no any kind of conditions are asked on generalized functions in the neighbourhood of their
singularities.
\\
In fact, it is precisely due to the lack of the mentioned type of constraints that the nowhere
dense algebras, and more generally, a large class of the space-time foam algebras have a
flabby sheaf structure, see next subsection.

\subsection{The Flabby Sheaf Property}\label{subsec2.3}

Regarding \textit{singularities}, it is particularly important to point out the following, see
Rosinger [7, pp. 173-187]. There exists a close connection between the amount of singularities
given spaces of generalized functions can handle, and whether the respective spaces form a
flabby sheaf.

\medskip

As a basic fact, let us recall that none of the spaces $\mathcal{C}^l(\mathbb{R} ^n)$, with
$0 \leq l \leq \infty$, form a flabby sheaf. And clearly, these spaces of usual functions of
various degrees of smoothness are not supposed to handle any of the more involved
singularities. If we now consider the \textit{smallest flabby sheaves} generated by these
spaces then, according to Oberguggenberger \& Rosinger [pp. 142-146], they are given, for $0
\leq l \leq \infty$, by

$$ \mathcal{C}^l_{nd}(\mathbb{R}^n) ~=~
             \{~ f : \mathbb{R}^n \longrightarrow \mathbb{R} ~~|~~
             \exists~ \Gamma ~\subset~ \mathbb{R}^n,~ \Gamma \textnormal{ closed, nowhere
             dense } :~ f \in \mathcal{C}^l(\mathbb{R}^n \setminus \Gamma) ~\} $$

And clearly, these spaces allow as singularities arbitrary closed nowhere dense subsets
$\Gamma ~\subset~ \mathbb{R}^n$, subsets which therefore can have arbitrary large Lebesgue
measures, Oxtoby. Furthermore, there is no any restriction on the functions in these spaces in
the neighbourhood of such singularities.

\medskip

Needless to say, the Schwartz distributions and the Colombeau algebras contain the spaces of
smooth functions $\mathcal{C}^l_{nd}(\mathbb{R}^n)$, with $0 \leq l < \leq \infty$. And in
view of the above, if the Schwartz distributions or the Colombeau algebras were flabby sheaves,
they would have to contain $ \mathcal{C}^l_{nd}(\mathbb{R}^n)$, thus they would have to be
able to deal with singularities given by arbitrary closed nowhere dense subsets, with no
conditions in their neighbourhoods on the respective generalized functions, and do so at least
in the case of usual functions.
\\
However, as is well known, such a class of singularities, plus the total lack of restrictions
in their neighbourhoods on the generalized functions, is not allowed for the Schwartz
distributions, or for that matter, in the Colombeau algebras. This failure, therefore, can be
seen as a consequence, among others, of such spaces of generalized functions not being flabby
sheaves.

\bigskip

Let us now turn to the space-time foam algebras. First, we define a large class of space-time
foam algebras $B_{L,~\mathcal{S}}(X)$ on arbitrary finite dimensional smooth manifolds $X$,
each of which will in the next Theorem prove to have a \textit{fine sheaf} structure. From
this proof it will also follow that the respective algebras are \textit{flabby sheaves} as
well, in case they satisfy a further rather natural condition. This large class of algebras
contains the nowhere dense algebras, and thus the result presented here is a significant
extension of the similar recent result in Mallios \& Rosinger [1, Lemma 2], see also Mallios
\& Rosinger [2,3].

\bigskip

Given a family $\mathcal{S}$ of singularity sets $\Sigma \subset X$ for which the conditions
(\ref{label5}), (\ref{label6}) hold, we call that family \textit{locally finitely additive},
if and only if it satisfies also the condition :

\medskip

For any sequence of singularity sets $\Sigma_l \in \mathcal{S}$, with $l \in \mathbb{N}$, if
we take $\Sigma = \bigcup_{l \in \mathbb{N}}\,\Sigma_l$, then for every nonvoid open subset $V
\subseteq X$, we have $\Sigma \cap V \in \mathcal{S} |_{V}$, whenever

\begin{equation}\label{label2.5}
      \begin{array}{l}
            \forall  \quad x \in V ~: \\ \\
            \exists \quad \Delta ~\subseteq~ V,~ \Delta
                            \textnormal{ neighbourhood of } x ~: \\ \\
            \quad \quad \{~ l \in \mathbb{N} ~|~ \Sigma_l \cap \Delta \neq \phi ~\}
                               \textnormal{ is a finite set of indices}
       \end{array}
\end{equation}

It is easy to verify that, see (\ref{label2}), (\ref{label3}), the families of singularities
$\mathcal{S}_{nd}(X)$ and $\mathcal{S}_{\textit{Baire I}}(X)$ are both locally finitely
additive.
\\
Indeed, $\mathcal{S}_{\textit{Baire I}}(X)$ is trivially so, since any countable union of
first Baire category sets is still of first Baire category. As far as $\mathcal{S}_{nd}(X)$ is
concerned, it suffices to note two facts. First, a subset of a topological space is closed and
nowhere dense, if and only if it so in the neighbourhood of every point. Second, a finite
union of closed nowhere dense sets is again closed and nowhere dense.

\bigskip

Let us also note that $\mathcal{S}_{nd}(X)$ and $\mathcal{S}_{\textit{Baire I}}(X)$ are among
those classes of singularities $\mathcal{S}$ which for every nonvoid open subset $V \subseteq
X$ satisfy the condition

\begin{equation}\label{label2.6}
\mathcal{S} | _{V} ~\subseteq~ \mathcal{S}
\end{equation}

We shall use the concept of \textit{sheaf} as defined by its \textit{sections}, see Bredon,
Oberguggenberger \& Rosinger, Mallios \& Rosinger [1,2], specifically we shall deal with
\textit{restriction} mappings to nonvoid open subsets $V \subseteq X$, or in particular, given
by the charts $V \in \mathcal{V}$.

\bigskip

Let us assume given a family $\mathcal{S}$ of singularities which satisfies the conditions
(\ref{label5}), (\ref{label6}). Then it is clear that for every nonvoid open subset $V
\subseteq X$ the restriction $\mathcal{S} |_{V}$ of $\mathcal{S}$ to $V$, see (\ref{label2.4}),
will also satisfy (\ref{label5}), (\ref{label6}), this time however on $V$.

\bigskip

Let now be given any right directed partial order $L=(\Lambda, \leq)$. Then the restriction to
all the nonvoid open subsets $V \subseteq X$ of the space-time foam algebra
$B_{L,~\mathcal{S}}(X)$ is the family of space-time foam algebras

\begin{equation}\label{label2.7}
\mathcal{B}_{L,~\mathcal{S}, X} ~=~ (~ B_{L,~\mathcal{S}|_{V}}(V) ~~|~~ V
~\subseteq~ X,~ V \textnormal{ nonvoid open} ~)
\end{equation}

a relation which follows easily, if we take into account (\ref{label2.4}), and the fact that

\begin{equation}\label{label2.8}
B_{L,~\mathcal{S}}(X)|_{V} ~=~ B_{L,~\mathcal{S}|_{V}}(X)
\end{equation}

which is a direct consequence of (\ref{label12}), (\ref{label10}), as well as of the obvious
relation, see (\ref{label8}),

\begin{equation}\label{label2.9}
\mathcal{J}_{L,~\Sigma}(X)|_{V} ~=~ \mathcal{J}_{L,~\Sigma\cap V}(V),~~
\textnormal{ for } \Sigma\subseteq X
\end{equation}

We can also note that in case $\Sigma \cap V = \phi$, the ideal $\mathcal{J}_{L,~\phi}(V)$,
and thus the algebra $B_{L,~\phi}(V)$, are still well defined, as long as $V$ is open and
nonvoid, see (\ref{label8}), (\ref{label10}), (\ref{label12}).
\\

\textbf{Theorem}

\medskip

Given on a finite dimensional smooth manifold $X$ any family of singularities $\mathcal{S}$
which is locally finitely additive.
\\
Then for every right directed partial order $L=(\Lambda, \leq)$, the family of space-time foam
algebras, see (\ref{label2.7})

\begin{equation}\label{label2.10}
\mathcal{B}_{L,~\mathcal{S},~X} ~=~ (B_{L,~\mathcal{S}|_{V}}(V) ~~|~~ V
\subseteq X,~ V \textnormal{ nonvoid open} ~)
\end{equation}

is a fine sheaf on $X$.

\bigskip

If in addition $\mathcal{S}$ has the property

\begin{equation}\label{label2.11}
     \begin{array}{l}
          \forall \quad \Sigma \in \mathcal{S},~ V \subseteq X,~ V
                    \textnormal{ nonvoid open},~ \Gamma \in \mathcal{S}_{nd}(V) ~: \\ \\
          \quad \quad (\Sigma \cap V) \cup \Gamma \in \mathcal{S}|_{V}
     \end{array}
\end{equation}

and furthermore, $\mathbb{N}$ is, see subsection \ref{subsec1.10}., cofinal in $\Lambda$, then
$\mathcal{B}_{L,~\mathcal{S},~X}$ in (\ref{label2.10}) is as well a flabby sheaf on $X$.

\bigskip

\textbf{Note 1}

\bigskip

The classes of singularities $\mathcal{S}_{nd}(X)$ and $\mathcal{S}_{\textit{Baire I}}(X)$
satisfy the conditions in the above Theorem.

\bigskip

\textbf{Note 2}

\bigskip

If we consider $\mathbb{N} \times \mathbb{N}$ with the partial order in subsection
\ref{subsec1.9}., and we embed $\mathbb{N}$ into $\mathbb{N} \times \mathbb{N}$ through the
diagonal mapping $\mathbb{N} \ni l \mapsto (l,l) \in \mathbb{N} \times \mathbb{N}$, then
$\mathbb{N}$ will be \textit{cofinal} in $\mathbb{N} \times \mathbb{N}$. Thus in view of the
above Theorem and Note 1, it follows that

$$\mathcal{B}_{\mathbb{N} \times \mathbb{N},~\mathcal{S}_{\textit{Baire I}}(X),~X} ~=~
     (~ B_{L,~\mathcal{S}_{\textit{Baire I}}(X) |_{V}}(V) ~~|~~ V
     ~\subseteq~ X,~ V \textnormal{ nonvoid open} ~)$$

is a \textit{fine flabby sheaf}.
\\
This result is nontrivial since $\mathcal{S}_{\textit{Baire I}}(X)$ contains lots of
singularity sets $\Sigma \subseteq X$, which are both \textit{dense} in $X$ and
\textit{uncountable}. In particular, this result is a significant strengthening of an earlier
similar result in Mallios \& Rosinger [1, Lemma 2] , where it was only given in the case of
the family of singularities $\mathcal{S}_{nd}(X)$.

\bigskip

\textbf{Note 3}

\bigskip

In the context of \textit{flabbiness} of spaces of functions or generalized functions, the
presence of $\mathcal{S}_{nd}(X)$ in condition (\ref{label2.11}) appears to be quite natural.
For instance, as seen in Oberguggenberger \& Rosinger [Remark 7.5, pp. 142-146], the class
$\mathcal{S}_{nd}(X)$ of closed nowhere dense singularities appears when one constructs the
\textit{smallest flabby sheaf} which contains $\mathcal{C}^{\infty}(X)$, for a nonvoid open
subset $X \subseteq \mathbb{R}^n$. The same happens when constructing the smallest flabby
sheaf containing $\mathcal{C}^0(X)$.

\bigskip

\textbf{Note 4}

\bigskip

One of the first important applications of the Theorem above is in Mallios \& Rosinger [2,3],
where the Abstract Differential Geometry theory Mallios [1,2] is illustrated with the most
\textit{singular structure coefficients} so far, namely, coefficients which can have
\textit{dense} singularity sets, and in particular, sets in $\mathcal{S}_{\textit{Baire I}}
(X)$. Earlier, in Mallios \& Rosinger [1], a similar illustration was given with structure
coefficients having singularity sets which could be closed and nowhere dense.

\medskip

These structure coefficients are in their most general case given by the space-time foam
differential algebras of generalized functions, which now replace the usual smooth functions
used in standard Differential Geometry.

\medskip

The important point to note with respect to the Abstract Differential Geometry in Mallios
[1,2] is that it is convenient for the structure coefficients - which as noted, in that
theory replace the usual smooth functions - to form sheaves of algebras that are both
\textit{fine} and \textit{flabby}. And it is precisely here that the Theorem above comes in,
by providing large classes of such algebras of generalized functions, including those
whose elements can have dense singularity sets, among them, as mentioned, in
$\mathcal{S}_{\textit{Baire I}}(X)$.

\bigskip

\textbf{Proof}

\bigskip

It is easy to see that the restrictions in (\ref{label2.10}) will satisfy the required sheaf
conditions related to restriction to open subsets.
\\
We can thus turn to check whether (\ref{label2.10}) satisfies the first of the two conditions
on sheaves, related to open covers. For that, here and in the sequel of the proof, we shall
for convenience, and without loss of generality, assume that our atlas $\mathcal{V}$ contains
all nonvoid open subsets in $X$.
\\
Let us therefore take $V = \bigcup_{i\in I} V_i$, where $V_i \subseteq X$, with $i \in I$, are
nonvoid open, while $I$ is any index set. Given now two generalized functions $T, T' \in
B_{L,~\mathcal{S} |_{V}}(V)$, let us assume that for $i\in I$ , we have

\begin{equation}\label{label2.12}
T |_{V_i} ~=~ T' |_{V_i}
\end{equation}

We then prove that

\begin{equation}\label{label2.13}
T ~=~ T'
\end{equation}

Indeed, let us note the relations

\begin{equation}\label{label2.14}
T ~=~ t + \mathcal{J}_{L,~S|_{V}}(V), \quad T' = t' +
\mathcal{J}_{L,~S|_{V}}(V)
\end{equation}

with suitably chosen $t,t' \in (\mathcal{C}^{\infty}(V))^{\Lambda}$, which follow from
(\ref{label12}). Then (\ref{label2.12}) implies for $i \in I$

\begin{equation}\label{label2.15}
(t'-t)|_{V_i} ~=~ w_i ~=~ (~ w_{i, \, \lambda} ~|~ \lambda \in \Lambda ~) \in
\mathcal{J}_{L,~\mathcal{S}|_{V_i}}(V_i)
\end{equation}

Now for $i \in I$, we define the product mapping

\begin{equation}\label{label2.16}
\mathcal{C}^{\infty}_{V_i}(V) \times \mathcal{C}^{\infty}(V_i) \ni
(\alpha, \psi)\longrightarrow \alpha \psi \in
\mathcal{C}^{\infty}(V)
\end{equation}

where $\mathcal{C}^{\infty}_{V_i}(V)$ denotes the set of all smooth functions in
$\mathcal{C}^{\infty}(V)$ whose support is in $V_i$, while the product $\alpha \psi$ is
defined by

\begin{displaymath}
(\alpha \psi)(x) = \left| \begin{array}{ll}
                                \alpha(x)\psi(x) & \textnormal{ if $x\in V_i$} \\ \\
                                0 & \textnormal{ if $x\in V\setminus V_i$}\\
                          \end{array} \right.
\end{displaymath}

At this point, we consider a smooth partition of unity $(~ \alpha_l ~|~ l\in \mathbb{N} ~)$ on
$V$, such that, see de Rham, each $\alpha_l$ has a compact support contained in one of the
$V_i$, and in addition, every point of $V$ has a neighbourhood which intersects only a finite
number of the supports of the various $\alpha_l$. In this way we obtain a mapping

\begin{equation}\label{label2.17}
N \ni l ~\longrightarrow~ i(l)\in I, \textnormal{ with  supp $\alpha_l\subseteq V_{i(l)}$}
\end{equation}

Extending now termwise the above product to sequences of smooth functions indexed by $\lambda
\in \Lambda$, we can define the sequence of smooth functions, see (\ref{label2.15})

\begin{equation}\label{label2.18}
w ~=~ (~ w_{\lambda}~|~ \lambda \in \Lambda ~) ~=~ \sum_{l\in \mathbb{N}}\alpha_l
\, w_{i(l)} \in (\mathcal{C}^{\infty}(V))^{\Lambda}
\end{equation}

and then show that

\begin{equation}\label{label2.19}
w\in \mathcal{J}_{L,~\mathcal{S}|_{V}}(V)
\end{equation}

Once we have (\ref{label2.19}), we recall (\ref{label2.15}) - (\ref{label2.17}) and the fact
that $(~ \alpha_l ~|~ l \in \mathbb{N} ~)$ is a smooth partition of unity on $V$, and we
obtain

$$t'-t ~=~ (~ \sum_{l\in \mathbb{N}}\alpha_l ~)(t'-t)~=~ \sum_{l\in
\mathbb{N}}\alpha_l (t'-t) ~=~ \sum_{l\in \mathbb{N}}\alpha_l (t'-t)|_{V_{i(l)}} ~=~
\sum_{l\in \mathbb{N}}\alpha_l w_{i(l)} ~=~ w$$

which in view of (\ref{label2.14}) will indeed yield (\ref{label2.13}).

\medskip

In order to obtain (\ref{label2.19}), and in view of (\ref{label10}), (\ref{label2.4}), it
suffices to find $\Sigma \in \mathcal{S}$, for which we have

\begin{equation}\label{label2.20}
w \in \mathcal{J}_{L,~\Sigma \cap V}(V)
\end{equation}

Let us therefore recall (\ref{label2.15}), and note that together with (\ref{label2.16}),
(\ref{label2.17}) and (\ref{label8}), it results in

\begin{equation}\label{label2.21}
\alpha_l w_{i(l)}\in \mathcal{J}_{L,~\mathcal{S}|_V}(V), \quad
\textnormal{ for $l\in \mathbb{N}$}
\end{equation}

thus (\ref{label10}) gives a sequence of sets $\Sigma_l^V \in \mathcal{S}|_V$, with $l\in
\mathbb{N}$, such that

\begin{equation}\label{label2.22}
\alpha_l w_{i(l)} \in \mathcal{J}_{L,~\Sigma_l^V}(V), \quad \textnormal{ for $l\in \mathbb{N}$}
\end{equation}

and due to (\ref{label2.21}), (\ref{label2.17}), (\ref{label8}), we can further assume about
$\Sigma_l^V$ that

\begin{equation}\label{label2.23}
\Sigma_l^V ~\subseteq~ \textnormal{supp }\alpha_l, \quad \textnormal{ for $l\in \mathbb{N}$}
\end{equation}

since we can always replace in (\ref{label2.22}) the initial $\Sigma_l^V$ with $\Sigma_l^V
\cap \textnormal{supp $\alpha_l$}$.

\medskip

However, for $l\in \mathbb{N}$ we have $\Sigma_l^V=\Sigma_l \cap V$, with suitable $\Sigma_l
\in \mathcal{S}$. Then by taking

\begin{equation}\label{label2.24}
\Sigma ~=~ \bigcup_{l\in \mathbb{N}} \Sigma_l
\end{equation}

and recalling (\ref{label2.23}), we obtain for $l\in \mathbb{N},~ x \in V,~ \Delta \subseteq
V,~ \Delta$, neighbourhood of $x$, the relations

$$\Sigma_l \cap \Delta ~=~ ( \Sigma_l \cap V )\cap \Delta ~=~ \Sigma_l^V \cap \Delta
        ~\subseteq~ ( \textnormal{supp } \alpha_l ) \cap \Delta$$

It follows therefore from the assumed property of the supports of the partition of unity
$(~ \alpha_l ~|~ l \in \mathbb{N} ~)$, that for the given $V$, the sequence of singularity
sets $\Sigma_l \in \mathcal{S}$, with $l\in \mathbb{N}$, satisfies condition (\ref{label2.5}).
Thus we have $\Sigma\cap V \in \mathcal{S}|_V$.

\medskip

But (\ref{label2.24}) yields

$$\Sigma \cap V ~=~ \bigcup_{l\in \mathbb{N}}(\Sigma_l \cap V)~=~
                                \bigcup_{l\in \mathbb{N}}\Sigma_l^V$$

and thus (\ref{label2.18}), (\ref{label2.22}) and (\ref{label8}), will give (\ref{label2.20}),
and in this way, the proof of (\ref{label2.18}) is completed.

\bigskip

As a last step in order to show that (\ref{label2.10}) is a sheaf, let $T_i \in
B_{L,~\mathcal{S}|_{V_i}}(V_i)$, with $i\in I$, be such that

\begin{equation}\label{label2.25}
T_i |_{V_i \cap V_j} ~=~ T_j |_{V_i \cap V_j}
\end{equation}

for all $i,j\in I$, for which $V_i \cap V_j \neq \phi$.  Then we show that

\begin{equation}\label{label2.26}
\exists~~~ T \in B_{L,~\mathcal{S} |_V}(V) ~: ~~\forall~~~ i \in I ~:  ~~~~T |_{V_i}~=~ T_i
\end{equation}

Indeed, (\ref{label12}) gives the representations

\begin{equation}\label{label2.27}
T_i ~=~ t_i + \mathcal{J}_{L,~\mathcal{S}|_{V_i}}(V_i), \quad \textnormal{ for $i\in I$}
\end{equation}

where $t_i \in (\mathcal{C}^{\infty}(V_i))^{\Lambda}$.  But then (\ref{label2.25}) results in

\begin{equation}\label{label2.28}
(t_i - t_j)|_{V_i \cap V_j} ~=~ w_{i, \, j}\in \mathcal{J}_{L,~\mathcal{S}|_{V_i
\cap V_j}}(V_i \cap V_j)
\end{equation}

for all $i,j \in I$ such that $V_i \cap V_j \neq \phi$.

\medskip

Let us now take any fixed $i\in I$. Given $l \in \mathbb{N}$ such that $V_i \cap V_{i(l)} \neq
\phi$, the relation (\ref{label2.28}) yields

$$t_{i(l)} ~=~ t_i + w_{i(l),\, i} \textnormal{ on $V_i \cap V_{i(l)}$}$$

thus (\ref{label2.16}), (\ref{label2.17}) lead to

$$\alpha_l t_{i(l)} ~=~ \alpha_l t_i + \alpha_l w_{i(l)\, i} \textnormal{ on $V_i$}$$

But then

$$ \sum_{l \in \mathbb{N}} \alpha_l t_{i(l)}~=~
           \left( \sum_{l \in \mathbb{N}}\alpha_l \right)t_i +
             \sum_{l\in \mathbb{N}} \alpha_l \, w_{i(l), \, i}~~~~~~ \textnormal{ on $V_i$} $$

or

\begin{equation}\label{label2.29}
\sum_{l \in \mathbb{N}} \alpha_l t_{i(l)}~=~
t_i + \sum_{l \in \mathbb{N}} \alpha_l w_{i(l),\, i}~~~~~~ \textnormal{ on $V_i$}
\end{equation}

On the other hand, the relation

\begin{equation}\label{label2.30}
\left(\sum_{l\in \mathbb{N}} \alpha_l w_{i(l),\, i}\right)\Bigg|_{V_i}\in
\mathcal{J}_{L,\, \mathcal{S}|_{V_i}}(V_i)
\end{equation}

follows by an argument similar with the one we used for obtaining (\ref{label2.19}), via
(\ref{label2.20}) - (\ref{label2.24}). In this way, if we define

\begin{equation}\label{label2.31}
T ~=~ t + \mathcal{J}_{L,\, \mathcal{S}|_V}(V) \in
B_{L,\, \mathcal{S}|_V}(V)
\end{equation}

where

$$ t ~=~ \sum_{l \in \mathbb{N}}\alpha_l t_{i(l)} $$

then (\ref{label2.29}) - (\ref{label2.31}) will give us (\ref{label2.26}), and the proof of
the fact that (\ref{label2.10}), is a sheaf is completed.

\bigskip

We turn now to proving that (\ref{label2.10}) is a \textit{fine} sheaf. This however follows
easily from (\ref{label16}), which as we have noted, implies that $1_V ~=~ u(1) +
\mathcal{J}_{L,\,\mathcal{S}|_V}(V)$ is the unit element in $B_{L,\,\mathcal{S}|_V}(V)$, thus
the partition of unity property together with (\ref{label8}) lead to

$$ 1_V ~=~ \sum_{l \in \mathbb{N}}(u(\alpha_l) + \mathcal{J}_{L,\,\mathcal{S}|_V}(V))
                     \in B_{L,\,\mathcal{S}|_V}(V) $$

At this point, we are left only with showing that (\ref{label2.10}) is a \textit{flabby} sheaf.

\medskip

Let $V'\subseteq V \subseteq X$ be nonvoid open subsets, and let $T' \in
B_{L,\,\mathcal{S}|_{V'}}(V')$. Let us denote by $\Sigma '$ the boundary of $V'$ in $V$. Then
clearly $\Sigma '$ is closed and nowhere dense in $V$, while $V \setminus (V' \cup \Sigma ')$
is open in $V$. Further, since $\Sigma '$ is closed in $V$, there exists, Kahn, $\sigma' \in
\mathcal{C}^{\infty}(V)$, such that $\Sigma' = \{~ x \in V ~|~ \sigma'(x)=0 ~\}$.
\\
We shall use now an auxiliary function $\eta \in \mathcal{C}^{\infty}(\mathbb{R})$ such that
$\eta = 1$ on $(-\infty,-1] \cup [1,\infty)$, while $\eta = 0$ on $[-1/2,1/2]$. And with its
help we can define the sequence of smooth functions $\beta_l \in \mathcal{C}^{\infty}(V)$,
with $l \in \mathbb{N}$, according to

\begin{equation}\label{label2.32}
\beta_l(x) ~=~ \left| \begin{array}{ll}
\eta((l+1)\sigma'(x)) & \textnormal{ if $x\in V'\cup \Sigma'$}\\\\
0 & \textnormal{if $x \in V \setminus (V'\cup \Sigma')$}\\
\end{array} \right.
\end{equation}

It is easy to check that

\begin{equation}\label{label2.33}
\textnormal{supp }\beta_l ~\subseteq~ V', \textnormal{ for $l \in \mathbb{N}$}
\end{equation}

and

\begin{equation}\label{label2.34}
\forall \, K \subset \subset V' \, :\, \exists \, l \in \mathbb{N}\,:\, \forall \, k \in
\mathbb{N},~ k\geq l \, : \, \beta_k =1~~~ \textnormal{ on } K
\end{equation}

Let us now assume that $T'$ has the representation, see (\ref{label12})

\begin{equation}\label{label2.35}
T' ~=~ t' + \mathcal{J}_{L,\, \mathcal{S}|_{V'}}(V')
\end{equation}

where $t'= ( t'_{\lambda} ~|~ \lambda \in \Lambda ) \in (\mathcal{C}^{\infty}(V'))^{\Lambda}$.
Then we define

\begin{equation}\label{label2.36}
T~=~ t + \mathcal{J}_{L,\, \mathcal{S}|_V}(V) \in B_{L,\, \mathcal{S}|_V}(V)
\end{equation}

where $t ~=~ ( t_\lambda ~|~ \lambda \in \Lambda )$ and, see (\ref{label2.16}),
(\ref{label2.33}), $t_\lambda ~=~ \beta_{l_\lambda}t'_{\lambda}$, for $\lambda \in \Lambda,~
l_{\lambda} \in \mathbb{N}$ and $\lambda \leq l_{\lambda}$, this last inequality being
possible, since we assumed that $\mathbb{N}$ is cofinal in $\Lambda$.

\bigskip

Before going further, we have to show that the definition of $T$ in (\ref{label2.36}) does not
depend on the choice of $t '$ in (\ref{label2.35}). Let us therefore assume that, instead of
the one in (\ref{label2.35}), we are given another representation

$$T' ~=~ t'^{*} + \mathcal{J}_{L,\, \mathcal{S}|_{V'}}(V')$$

with $t'^{*} = ( t'^{*}_{\lambda} ~|~ \lambda \in \Lambda ) \in (\mathcal{C}^{\infty}(V'))^
{\Lambda}$, then clearly

\begin{equation}\label{label2.37}
t'^{*} - t' ~=~ ( t'^{*}_{\lambda}-t'_{\lambda} ~|~ \lambda \in \Lambda ) \in
\mathcal{J}_{L,\, \mathcal{S}|_{V'}}(V')
\end{equation}

As above with $t$ in (\ref{label2.36}), let us now define $t^{*} = ( t^{*}_{\lambda} ~|~
\lambda \in \Lambda ) \in (\mathcal{C}^{\infty}(V))^{\Lambda}$ by $t^{*}_{\lambda} =
\beta_{l_{\lambda}} t'^{*}_{\lambda}$, for $\lambda \in \Lambda$. We show then that $T$ in
(\ref{label2.36}) has also the representation

$$ T ~=~ t^{*} + \mathcal{J}_{L,\, \mathcal{S}|_V}(V) \in B_{L,\, \mathcal{S}|_V}(V)$$

or that, equivalently

\begin{equation}\label{label2.38}
t^{*} - t \in \mathcal{J}_{L,\, \mathcal{S}|_V}(V)
\end{equation}

Indeed, from (\ref{label2.37}), (\ref{label10}), (\ref{label2.4}) we obtain $\Sigma \in
\mathcal{S}$ such that

$$ t'^{*} - t' \in \mathcal{J}_{L,\, \Sigma\cap V'}(V') $$

and then (\ref{label8}) gives

\begin{equation}\label{label2.39}
\begin{array}{l}
\forall~~ \, x \in V' \setminus(\Sigma \cap V') ~:\\ \\
\exists~~ \, \lambda \in \Lambda ~:\\ \\
\forall~~ \, \mu \in \Lambda,~ \mu \geq \lambda ~:\\ \\
\forall~~ \, p \in \mathbb{N}^p ~:\\ \\
\quad D^p(t'^{*}_{\lambda}-t'_{\lambda})(x)=0
\end{array}
\end{equation}

But in view of (\ref{label2.11}), it follows that $(\Sigma \cap V) \cup \Sigma' \in
\mathcal{S}|_V$, since clearly $\Sigma' \in \mathcal{S}_{nd}(V)$.  Also, we obviously have
$t^{*}_{\lambda} - t_{\lambda} = \beta_{l_{\lambda}}(t'^{*}_{\lambda} - t'_{\lambda})$, with
$\lambda \in \Lambda$. And then (\ref{label2.39}), (\ref{label2.32}) and (\ref{label8}) will
directly lead to

$$ t^{*} - t \in \mathcal{J}_{L,\, (\Sigma \cap V)\cup \Sigma '}(V)
                           ~\subseteq~ \mathcal{J}_{L,\, \mathcal{S}|_V}(V) $$

and the proof of (\ref{label2.38}) is completed.

\bigskip

At last, it follows easily from (\ref{label2.32}) - (\ref{label2.36}) that

$$ T|_{V'} ~=~ T' $$

since a direct computation using also (\ref{label8}) , gives

$$ t' - t|_{V'} \in \mathcal{J}_{L,\, \phi}(V') $$

In this way the flabbiness of (\ref{label2.10}) is proved.

\section{Connections with Distributions}

Let us indicate in short the way the multi-foam algebras can be related to the Schwartz
distributions. For that, let us recall in some detail the mentioned wide ranging purely
algebraic characterization of all those differential algebras of generalized functions in
which one can embed linearly the Schwartz distributions, a characterization which, as also
mentioned, contains the Colombeau algebras as a particular case, see Rosinger [4, pp. 75-88],
Rosinger [5, pp. 306-315], or Rosinger [6, pp. 234-244].

\bigskip

According to the mentioned characterisation, in the case of the multi-foam algebras
$B_{L,\, \mathcal{S}}(X)$, for instance, the \textit{necessary and sufficient} condition for
the existence of such a linear embedding, namely

\begin{equation}\label{label3.1}
\mathcal{D}'(X) ~\subset~ B_{L,\, \mathcal{S}}(X)
\end{equation}

is precisely the \textit{off-diagonality} condition (\ref{label15}), which as we have seen,
does indeed hold. Furthermore, the linear embedding (\ref{label3.1}) will preserve the
differential structure of $\mathcal{C}^{\infty}(X)$.

\section{Final Remarks}

\textbf{Remark 1.}

\bigskip

It is important to note that, just like in Mallios \& Rosinger[1] , where the nowhere dense
differential algebras of generalized functions were used, or for that matter, in Rosinger
[1-11], Colombeau, Oberguggenberger, Grosser et.al., where other differential algebras of
generalized functions appeared as well, so in this paper, where the space-time foam,
multi-foam or foam differential algebras of generalized functions are employed, there is again
\textit{no need} for any topological algebra structure on these algebras.
\\
This \textit{lack of need} of any kind of topological algebra structures on differential
algebras of generalized functions may appear somewhat strange to those accustomed with the
rather involved functional analytic methods. However, this lack of need is precisely one of
the \textit{strong} points of the nonlinear algebraic method of generalized functions
initiated and developed in the mentioned references, see in this regard 46F30 in the AMS
Subject Classification 2000 at www.ams.org/msc/46Fxx.html. In fact, such an approach which
does not make use of topological structures on the differential algebras of generalized
functions may be seen as placing the subject in what is called Algebraic Analysis, an approach
to Analysis which has a long, successful and respectable tradition, see Synowiec.
\\
One of the reasons for the lack of need for any topological algebra structure on the algebras
of generalized functions under consideration is the following. It is becoming more and more
clear that the classical Hausdorff-Kuratowski-Bourbaki topological concept is not suited to
the mentioned algebras of generalized functions. Indeed, these algebras prove to contain
plenty of \textit{nonstandard} type of elements, that is, elements which in a certain sense
are infinitely small, or on the contrary, infinitely large. And in such a case, just like in
the much simpler case of nonstandard reals $^*\mathbb{R}$, any topology which would be
Hausdorff on the whole of the algebras of generalized functions, would by necessity become
discrete, therefore trivial, when restricted to usual, standard smooth functions, see for
details Biagioni.

\bigskip

Here, in order further to clarify the issue of the possible limitations of the usual
Hausdorff-Kuratowski-Bourbaki concept of topology, let us point out the following. Fundamental
results from Measure and Integration Theory, predating the mentioned concept of topology, yet
having a clear topological nature, have never been given a suitable formulation within that
Hausdorff-Kuratowski-Bourbaki concept. Indeed, such is the case, among others, with the
Lebesgue dominated convergence theorem, the Lusin theorem on the approximation of measurable
functions by continuous ones, and the Egorov theorem on the relation between point-wise and
uniform convergence of sequences of measurable functions.
\\
Similar limitations of the Hausdorff-Kuratowski-Bourbaki concept of topology appeared in the
early 1950s, when attempts were made to turn the convolution of Schwartz distributions into an
operation simultaneously continuous in both its arguments. More generally, it is well known
that, given a locally convex topological vector space, if we consider the natural bilinear
form defined on its Cartesian product with its topological dual, then there will exist a
locally convex topology on this Cartesian product which will make the mentioned bilinear form
simultaneously continuous in both of its variables, if and only if our original locally convex
topology is in fact as particular, as being a normed space topology.
\\
It is also well known that in the theory of ordered spaces, in particular, ordered groups or
vector spaces, there are important concepts of convergence, completeness, etc., which have
never been given a suitable formulation in terms of the usual Hausdorff-Kuratowski-Bourbaki
concept of topology. In fact, as seen in Oberguggenberger \& Rosinger, powerful general
results can be obtained about the existence of generalized solutions for very large classes of
nonlinear PDEs, by using order structures and their Dedekind type order completions, without
any recourse to associated topologies. And the generalized solutions thus obtained have the
universal regularity property that they can be assimilated with usual measurable functions.
\\
Finally, it should be pointed out that, recently, differential calculus was given a new
foundation by using standard concepts in category theory, such as naturalness. This approach
also leads to topological type processes, among them the so called toponomes or
$\mathcal{C}$-spaces, which prove to be more general than the usual
Hausdorf-Kuratowski-Bourbaki concept of topology, see Nel, and the references cited there.

\bigskip

In this way, we can conclude that mathematics contains a variety of important
\textit{topological type processes} which, so far, could not be formulated in convenient terms
using the Hausdorff-Kuratowski-Bourbaki topological concept. And the differential algebras of
generalized functions, just as much as the far simpler nonstandard reals $^*\mathbb{R}$,
happen to exhibit such a class of topological type processes.
\\
Further details in this regard can be found in Rosinger \& Van der Walt.

\bigskip

On the other hand, the topological type processes on the nowhere dense differential algebras
of generalized functions, used in Mallios \& Rosinger [1], for instance, as well as on the
space-time foam or special space-time foam differential algebras of generalized functions
employed in this paper, can be given a suitable formulation, and correspondingly, treatment,
by noting that the mentioned algebras are in fact \textit{reduced powers}, see Lo\u{s},
or Bell \& Slomson, of $\mathcal{C}^{\infty}(X)$, and thus of $\mathcal{C}(X)$ as well. Let us
give some further details related to this claim in the case of the space-time foam algebras,
see also Mallios \& Rosinger [2]. The case of the nowhere dense algebras was treated in
Mallios \& Rosinger [1].
\\
Let us recall, for instance, the definition in (\ref{label12}) of the multi-foam algebras, and
note that it obviously leads to

\begin{equation}\label{label4.1}
\begin{array}{l}
B_{L,\, \mathcal{S}}(X) ~=~
(\mathcal{C}^{\infty}(X))^{\Lambda} / \mathcal{J}_{L,\, \mathcal{S}}(X) ~\subseteq~ \\ \\
\quad \subseteq~ (\mathcal{C}(X))^{\Lambda}/\mathcal{J}_{L,\, \mathcal{S}}(X)
~\subseteq~ \mathcal{C}(\Lambda \times X)/\mathcal{J}_{L,\, \mathcal{S}}(X)
\end{array}
\end{equation}

assuming in the last term that on $\Lambda$ we consider the discrete topology.

\bigskip

Now it is well known, Gillman \& Jerison, that the algebra structure of $\mathcal{C}(\Lambda
\times X)$ is connected to the topological structure of $\Lambda \times X$. However, this
connection is rather sophisticated, as essential aspects of it involve the Stone-\u{C}ech
compactification $\beta(\Lambda \times X)$ of $\Lambda \times X$. And in order to complicate
things, in general $\beta(\Lambda \times X) \neq \beta(\Lambda) \times \beta(X)$, not to
mention that $\beta(\Lambda)$ alone, even in the simplest nontrivial case of
$\Lambda = \mathbb{N}$, has a highly complex structure.
\\
It follows that a good deal of the discourse, and in particular, the topological type one, in
the space-time foam algebras may be captured by the topology of $\Lambda \times X$, and of
course, by the far more involved topology of $\beta(\Lambda \times X)$. Furthermore, the
differential properties of these algebras will, in view of (\ref{label17}) - (\ref{label19}),
be reducible termwise to the classical differentiation of sequences of smooth functions.

\bigskip

In short, in the case of the mentioned differential algebras of generalized functions, owing
to their structure of \textit{reduced powers}, one obtains a 'two-way street' along which, on
the one hand, the definitions and operations are applied to sequences of smooth functions, and
then reduced termwise to such functions, while on the other hand, all of that has to be done
in a way which will be compatible with the `reductions' of the `power' by the quotient
constructions in (\ref{label12}), or in other words, (\ref{label4.1}). By the way, such a
`two-way street' approach has ever since the 1950s been fundamental in the branch of
Mathematical Logic, called Model Theory, see Lo\u{s}.
\\
But in order not to become unduly overwhelmed by ideas of Model Theory, let us recall here
that the classical Cauchy-Bolzano construction of the real numbers $\mathbb{R}$ is also a
reduced power. Not to mention that a similar kind of reduced power construction - in fact, its
particular case called `ultra-power' - gives the nonstandard reals $^*\mathbb{R}$ as well.

\bigskip

\textbf{Remark 2.}

\bigskip

Lately, there has been a growing interest in \textit{noncommutative} studies, and in
particular algebras, see Connes. It is therefore appropriate to mention possible connections
between such noncommutative methods and the space-time foam differential algebras of
generalized functions in this paper.
\\
ln this regard, we recall that, as mentioned at the beginning, in case our constructions
starts with arbitrary Banach algebra valued, and not merely real or complex valued, functions
then the resulting space-time foam algebras can in genral become \textit{noncommutative}.
\\
On the other hand, the emergence of noncommutative studies need not at all mean the loss of
any interest in, and relevance of commutative structures. Indeed, in many problems the
commutative approach turns out to be both more effective and also, of course, much more
simple.

\bigskip

Finally, it is important to mention here that in the case of singularities of generalized
functions, that is, of singularities in a differential context, the approach in Connes falls
far short even of the long established linear theory of Schwartz distributions. Indeed, the
only differential type operation in Connes, see pp. 19-28, 287-291, is defined as the
commutator with a fixed operator. In this way, it is a distinctly particular, and in fact,
rather trivial derivation, even when considered within Banach algebras. The effect is that,
it can only to a small extent deal with the singularities which otherwise even the limited
linear Schwartz Distributions Theory can handle. And in fact, the approach in Connes can deal
even less with singularities on arbitrary closed nowhere dense sets, let alone, on the far
larger class of arbitrary dense sets, such as those in this paper.

%%%%%%%%%%%%%%%%%%%%%%%%%%%%%%%%%%%%%%%%%%%%%%%%%%%%%%%%%%%%%%%%%%%%%%%%%%%%%%%%%%%%%%%%%%%%%%%%%%%%%%%%%%%%%%

\bibliographystyle{amsplain}

% ----------------------------------------------------------------
\end{document}